\documentclass[12pt,a4paper]{article}
\usepackage[english]{babel}
\usepackage{amsmath}
\usepackage{graphicx}

\usepackage{bbm}
\usepackage[%
breaklinks,
colorlinks=true,
linkcolor=blue,
anchorcolor=blue,
citecolor=black,
filecolor=black,
menucolor=black,
urlcolor=blue]{hyperref}
\DeclareMathOperator{\sign}{sign}

\begin{document}

\title{ccc-Autoevolutes}
\author{ \href{https://www.math.uni-bonn.de/people/karcher/}{Hermann Karcher},
		\href{https://page.math.tu-berlin.de/~ekki}{Ekkehard-H. Tjaden}
            }

\date{Jan.\,28.\,2021}

\maketitle

\begin{abstract}
\noindent ccc-Autoevolutes are closed constant curvature space curves
which are their own evolutes.
A modified Frenet equation produces curves which are congruent to their
evolutes.
Using symmetries we construct closed curves by solving 2-parameter problems
numerically.
The images are made with the program 3D-XplorMath~\cite{3DXM}.

\noindent%
Classification: \href{http://msc2010.org/MSC-2010-server.html}{53A04}
\end{abstract}

\vglue-15pt
\begin{figure}[h]  
  \centering
 \hbox{
 \vbox{\hsize=0.48\hsize  \includegraphics[width=2.7in]{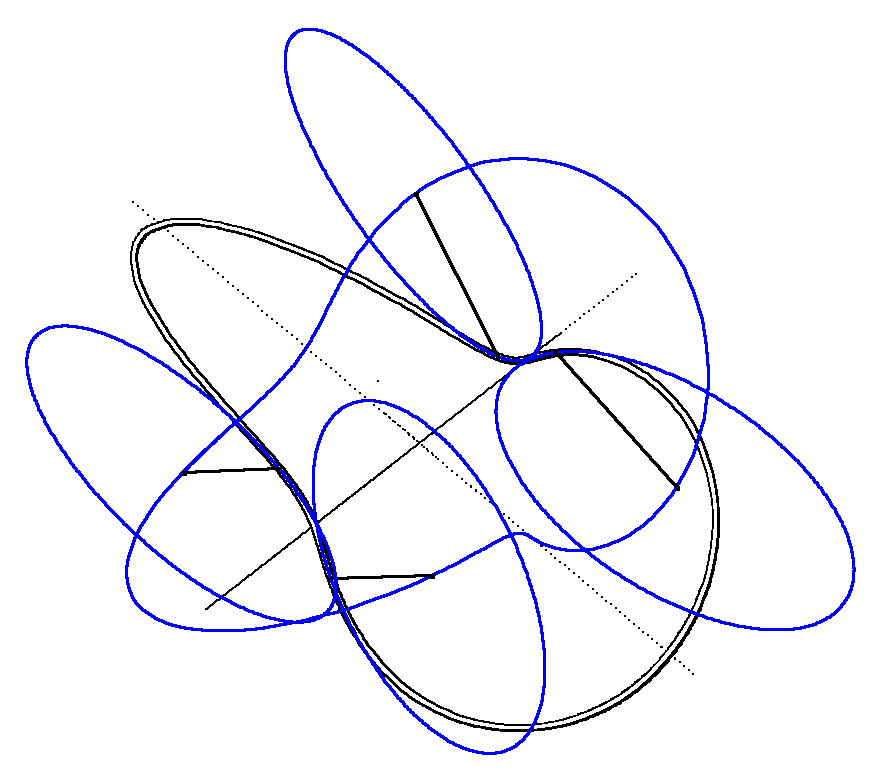} }
 \vbox{\hsize=0.452\hsize  \includegraphics[width=3in]{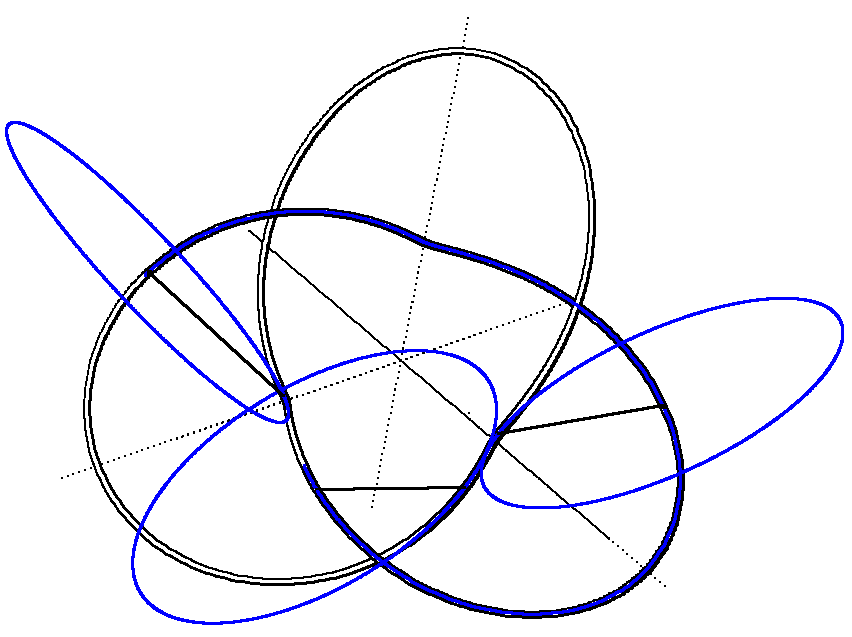} }
         }
  \caption{  \vbox{Left: constant curvature curve with congruent evolute \hfill \break
 			Right: 2--3--autoevolute, both with some osculating circles\vskip-17pt 
	      }           }
  \label{fig:ae0}
\end{figure}

Osculating circles moving along the curve can be seen in \cite{OscAE}.

\vskip10pt
This paper is a continuation of \cite{Karcher2020}. It started with the challenge from
Ekkehard-H. Tjaden to look for constant curvature curves which are congruent to their
evolutes. The details are joint work.

Assume that the curve $c$ is parametrized by arc length $s$,
has constant curvature $\kappa$ and its Frenet frame is
$\left\{ T(s), N(s), B(s) \right\}$.
Its evolute is
\begin{align}
  \label{eq:frenet1}
  \tilde{c}(s) &= c(s) + \frac{1}{\kappa}\,N(s)\\
  \tilde{c}\,'(s) &= \frac{\tau(s)}{\kappa}\,B(s)\,.
\end{align}
The velocity of the evolute is therefore
\begin{align}
  \label{eq:velo}
  \|\tilde{c}\,'(s)\| &= \frac{|\tau(s)|}{\kappa}
\end{align}
\vglue-10pt
with
\begin{align}
  \label{eq:frenet2}
  \tilde{T}(s) &= \sign(\tau)\,B(s)\\
  \label{eq:frenet2a}
  \tilde{N}(s) &= -N(s)\\
  \label{eq:frenet2b}
  \tilde{B}(s) &= \sign(\tau)\,T(s)\,.
\end{align}

The derivative of \eqref{eq:frenet2a} yields
\begin{align}
  -N'(s) &= \kappa\,T(s) -\tau(s)\,B(s)\\
  \tilde{N}'(s) &= \|\tilde{c}'(s)\|\, \left(-\tilde{\kappa}\,\tilde{T}(s) + \tilde{\tau}(s)\,\tilde{B}(s)\right)\\
    &= -\frac{\tau(s)\,\tilde{\kappa}}{\kappa}\,B(s) +\frac{\tau(s)\,\tilde{\tau}(s)}{\kappa}\,T(s)
\end{align}
By comparsion we get
\begin{align}
  \tilde{\kappa} = \kappa, \hskip10pt
  \tilde{\tau}(s)\, \tau(s) = \kappa^2
\end{align}

The two curves
  $c$, $\tilde{c}$
are geodesics on the canal surface which envelops the spheres of radius
  $\frac{1}{2\,\kappa}$
with midpoints on the curve 
{$ m(s) = (c(s)+\tilde{c}(s))/2$},
because their principle normals
  $N$, $\tilde{N}$
are orthogonal to the canal surface. This canal surface is, in the case of
the two examples above, a wobbly torus -- maybe that helps the visualization.

Constant curvature curves in general do not have congruent evolutes and even if
  $c$, $\tilde{c}$
are congruent, this is difficult to check in case
  $c$ 
is parametrized by arc length. One may
look for curves with non-constant velocity
  $v(t) = | \dot{c}(t)|$,
choose
  $\tau(t)$
in terms of
  $v(t)$
and try to arrange things so
that the first half of the curve is congruent to the evolute
of the second half and vice versa.
 As a consequence both halves of the middle curve $m(s)$
 are congruent.This is achieved by the
following version of the Frenet equations which Ekkehard
suggested. It assumes
  $t \in [0,2 \pi]$
and
  $v(t+\pi) = 1/v(t)$.
\begin{align}
  \label{eq:frenet3}
  \dot{T}(t) &= \kappa\,v(t)N(t)\\
  \dot{N}(t) &= -\kappa\,v(t)\,T(t) + \kappa/v(t)\,B(t)\\
  \dot{B}(t) &= -\kappa/v(t)\,N(t)\,.
\end{align}

For our computations we choose first
  $h(t) = a\, (\sin(t) + b_3\cdot\sin(3t))$
and then
  $v(t) = \sqrt{1+h(t)^2} -h(t)$
or also
  $v(t) = \exp(h(t))$.
We find it convenient not to restrict us to
  $\kappa = 1$, but to also vary the curvature.
Comparison with the standard Frenet equations shows
  $\tau(t) = \kappa/v(t)^2$.
Note that
  $v(t)$ and $\tau(t)$
are, relative to
  $t^* = \pi/2+n\,\pi$,
even functions. This implies that the principle normals at these
points are symmetry normals, saying that
  $180^\circ$
rotations around these normals map the curve onto itself.
Closed
examples can therefore be found by solving the following 2-parameter problem:\\
a) get the symmetry normals to lie in a plane -- they then automatically
  all intersect in one point.\\
b) get neighbouring symmetry lines to intersect with
  a rational angle -- preferably angles like
  $\pi/2$, $\pi/3$, $2\pi/3$\ldots .

This looks similar to the case of single constant curvature
curves. However in that case the constant Fourier term
of the torsion is a parameter that allows to solve problem
a) for all choices of other parameters. This allows to deal
with problem b) assuming that a) is already solved.
For the autoevolutes we have not found such a simplifying
parameter. By trial and error one has to roughly close the curve.
Only then does the automatic solution of
the 2-paramter problem work and produce a high accuracy
solution.

\begin{figure} [h]
  \centering \hbox{\hglue-3pt
 \vbox{ \hsize=0.33\hsize\includegraphics[width=1.8in]{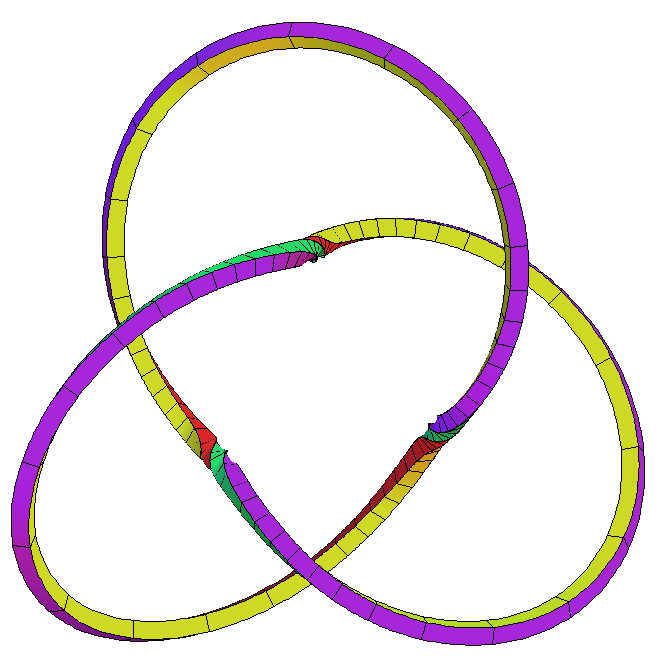}}
 \vbox{ \hsize=0.32\hsize\includegraphics[width=1.8in]{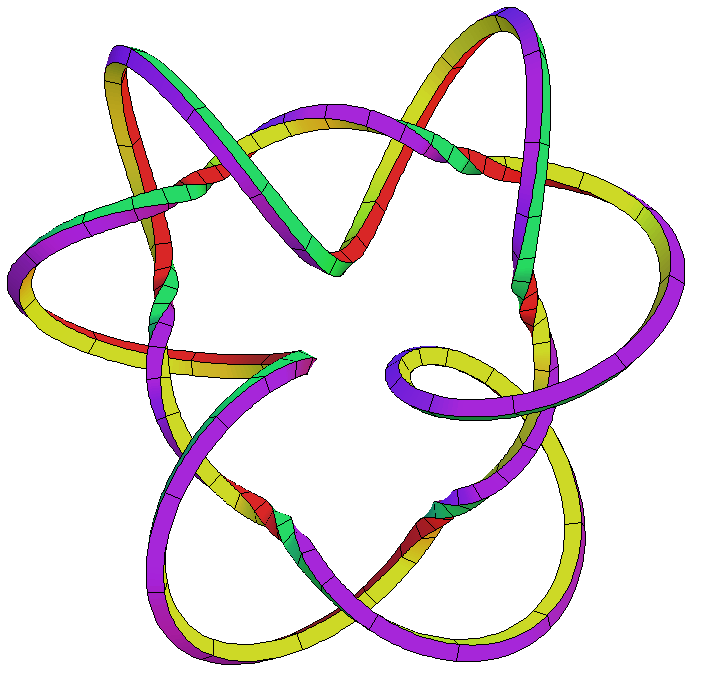}}
  \vbox{ \hsize=0.33\hsize\includegraphics[width=1.9in]{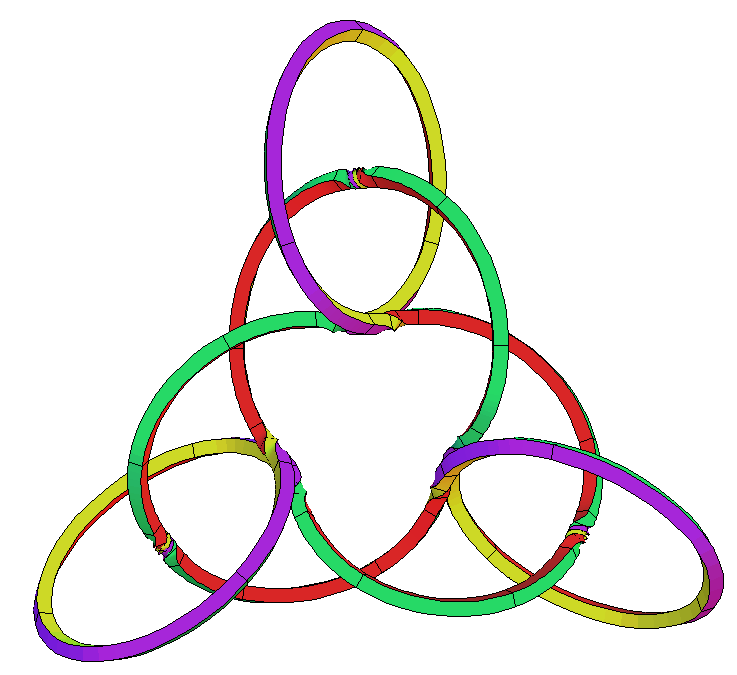}}
                            }
  \caption{
Note that the curves have large torsion – i.e. twist rapidly
and move slowly – opposite to the points with small
torsion, where the curves are fairly circular and their velocities are large.}
  \label{fig:ae1}
\end{figure}

Note that we are using only three parameters $\kappa, a, b_3$ to get these
examples. If we change $b_3$ in sufficiently small steps then we can adjust
$\kappa, a$ to deform the examples as closed curves. Of course, if we add
more (small) Fourier terms to the definition of $h(t)$, then we get higher
dimensional deformation families. By changing the parameters without 
adjustment we get to the vicinity of other examples, which then are caught
with our 2-dim equation solver. We found for odd $n$ and even $k$ 
autoevolutes which are $(n,k)$-torus knots. We found them with additional
inward loops as the second example above. The right example is still
recognizable with $n=5$, but for larger $n$ our examples are too complicated
to print. Here are some examples:
\vglue-8pt
\begin{figure} [h]
  \centering \hbox{\hglue-3pt
 \vbox{ \hsize=0.33\hsize\includegraphics[width=1.8in]{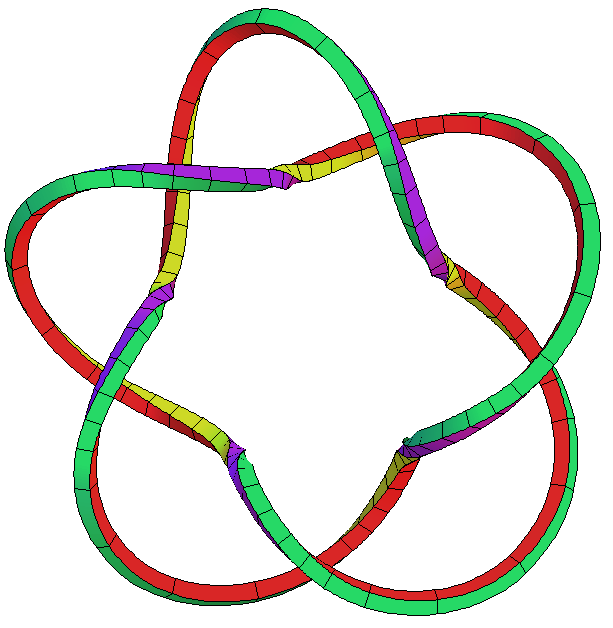}}
 \vbox{ \hsize=0.32\hsize\includegraphics[width=1.8in]{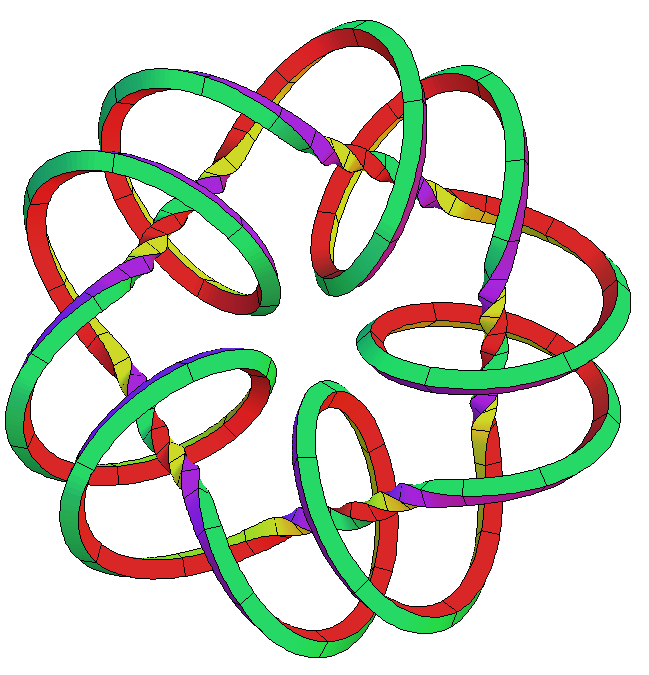}}
  \vbox{ \hsize=0.33\hsize\includegraphics[width=1.9in]{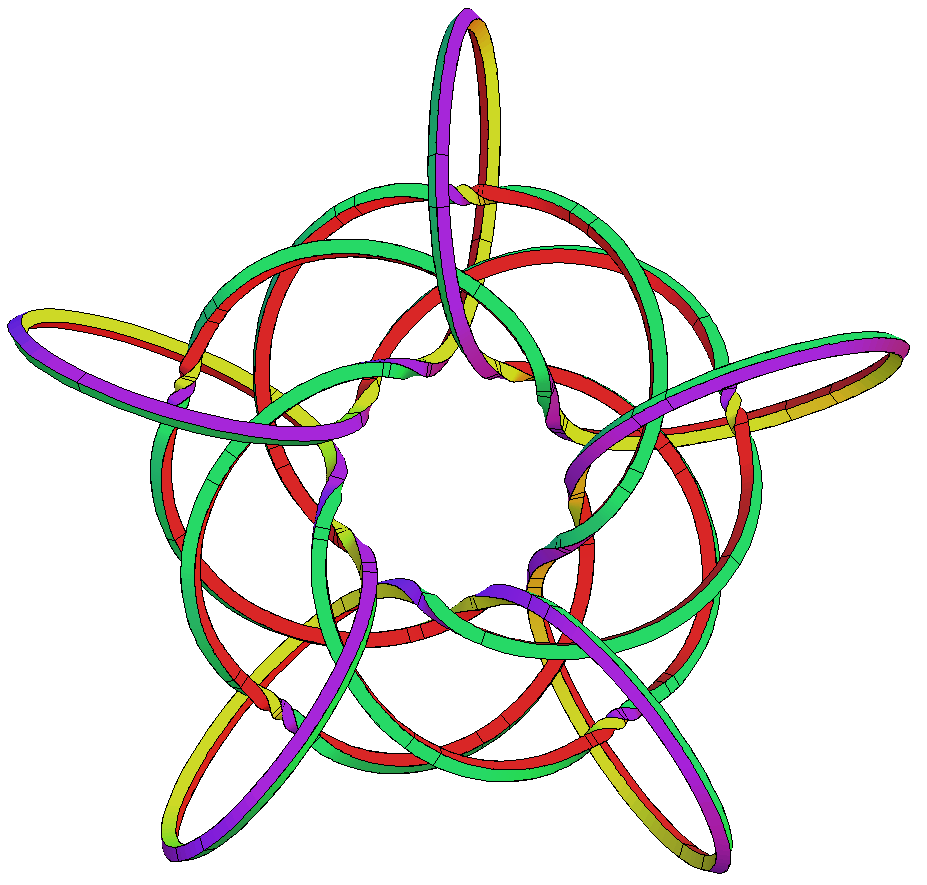}}
  			}
 \caption{2-5-knot, 2-5-knot with inward loops, 2-5-knot with outward loops. \\
 The canal surfaces on which these curves lie, look like a torus in the first case,
 but have heavy self-intersections in the other cases.}
\label{fig:ae2}                           
\end{figure} 
\vglue-16pt
\begin{figure} [h] 
   \centering \hbox{\hglue-3pt
 \vbox{ \hsize=0.33\hsize\includegraphics[width=1.9in]{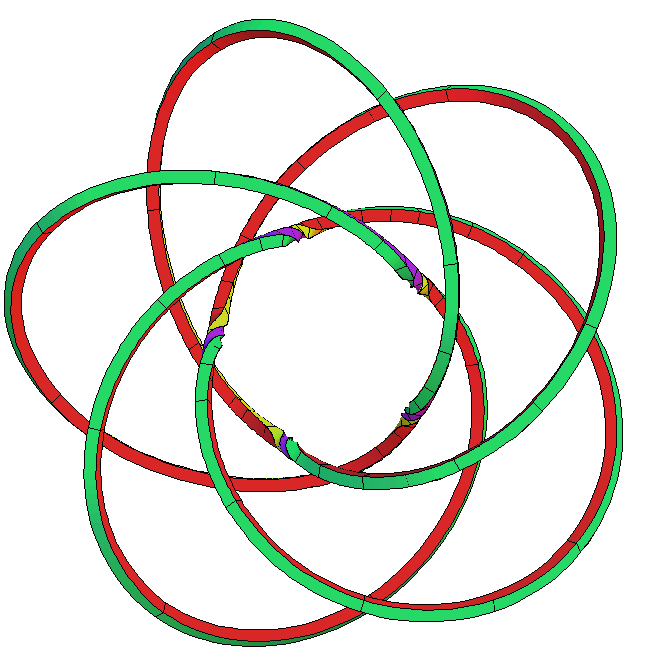}}
 \vbox{ \hsize=0.32\hsize\includegraphics[width=1.8in]{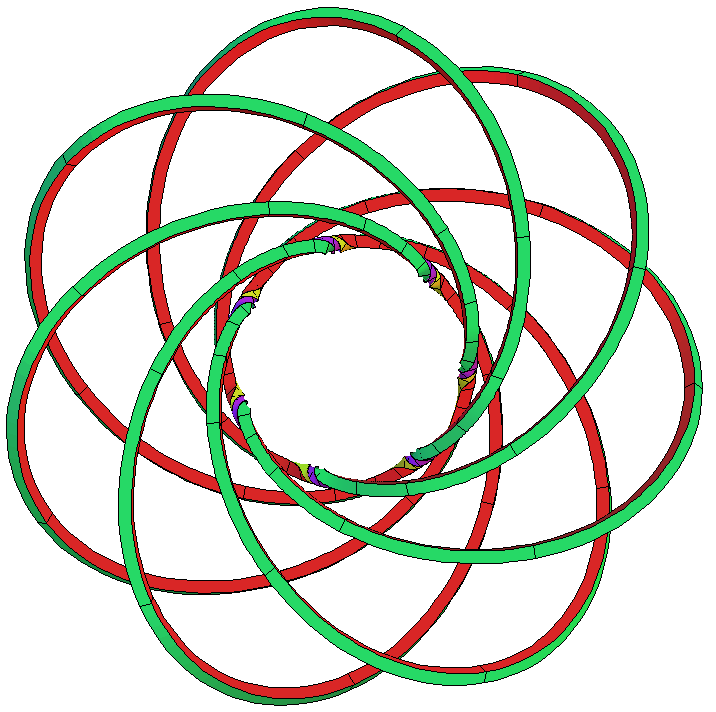}}
  \vbox{ \hsize=0.33\hsize\includegraphics[width=1.9in]{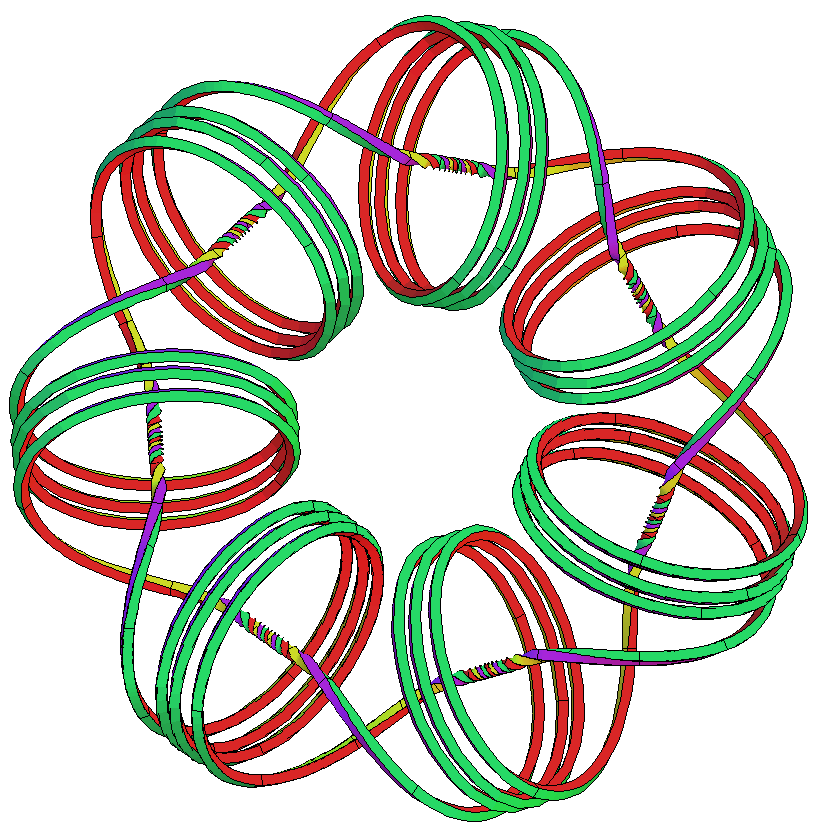}}
  			}
 \caption{4-5-knot, 6-7-knot,  2-7-knot with tripple inward loops. 
 Note the short pieces where the curve twists rapidly.}
\label{fig:ae3}                           
\end{figure}


\bibliographystyle{abbrv}
\bibliography{autoevolute3}

\begin{thebibliography}{10}

\bibitem{Karcher2020}
Hermann Karcher.
\newblock {\em Closed Constant Curvature Space Curves.\phantom{asdfghasd}}
\newblock \href{https://arxiv.org/abs/2004.10284}{https://arxiv.org/abs/2004.10284}

\bibitem{3DXM}
\newblock {\em Homepage of the program 3D-XplorMath}.
\newblock \href{http://3d-xplormath.org/}{http://3d-xplormath.org/}

\bibitem{OscAE}
{\em Osculating circles moving along a ccc-autoevolute.\phantom{asdfghasdasdfgasdfg}}
\newblock \href{https://www.math.uni-bonn.de/people/karcher/MyAEvideo.html}
               {https://www.math.uni-bonn.de/people/karcher/MyAEvideo.html} 

{\tt
Hermann Karcher $<$unm416@uni-bonn.de$>$\\
Ekkehard-H. Tjaden $<$tjaden@math.tu-berlin.de$>$
}

\end{thebibliography}
\goodbreak
\end{document}